\documentclass[12pt]{article}
\makeatletter
  
  \@addtoreset{equation}{section}
\makeatother
\newtheorem{proposition}{Proposition}
\newtheorem{theorem}{Theorem}
\newtheorem{lemma}{Lemma}
\newtheorem{remark}{Remark}
\newtheorem{corollary}{Corollary}
\newtheorem{definition}{Definition}
\usepackage{amsmath,amsfonts,latexsym}
\usepackage{color}
\date{ }

\begin{document}

\title{Stochastic optimal transport revisited
\thanks{Partially supported by JSPS KAKENHI Grant Numbers JP26400136 and 19K03548.
This is a final version which will appear in SN Partial Differential Equations and Applications.
}
}

\author{
Toshio Mikami\\
Dept. Math, Tsuda Univ.
}

\maketitle

\begin{abstract}

We prove the Duality Theorems for the stochastic optimal transportation problems with a convex cost function without a regularity assumption that is often supposed in the proof of the lower semicontinuity of an action integral.
In our new approach,
we prove that the stochastic optimal transportation problems with a convex cost function
are equivalent to a class of variational problems for the Fokker-Planck equation,
which lets us revisit them.
It is done by the so-called superposition principle and by an idea from the mather theory.
The superposition principle is the construction of a semimartingale from the Fokker-Planck equation
and can be considered a class of the so-called marginal problems that construct stochastic processes from given marginal distributions.
It was first considered in stochastic mechanics by E. Nelson, called Nelson's problem, and was proved by E. Carlen first.
The semimartingale is called the Nelson process, provided it is Markovian.
We also consider the Markov property of a minimizer of the stochastic optimal transportation problem with a nonconvex cost in a one-dimensional case.
In the proof, the superposition principle and the minimizer of an optimal transportation problem with a concave cost function play crucial roles.
Lastly, we prove the semiconcavity and the Lipschitz continuity of Schr\"odinger's problem that is a typical 
example of the stochastic optimal transportation problem.
\end{abstract}


\section{Introduction}
\label{intro}
The construction of a stochastic process from given marginal distributions is called a marginal problem.\\
\indent
Schr\"odinger's problem is the construction of a Markov diffusion process on $[0,1]$ from two endpoint marginal distributions at $t=0,1$
by solving a variational problem on the relative entropy.
We describe it briefly (see $V_S$ in (\ref{1.12}), (\ref{4.17}), and also \cite{1-Jamison1975,1-Leo2,Mikami}).
Let $\sigma$ and $\xi$ be, respectively, a $d\times d$ nondegenerate matrix-valued 
and an $\mathbb{R}^d$-valued function on $[0,1]\times \mathbb{R}^d$.
Suppose that the following stochastic differential equation has a weak solution $\{X(t)\}_{0\le t\le 1}$ 
with a positive transition probability density $p(s,x;t,y), 0\le s<t\le 1, x,y\in \mathbb{R}^d$:
\begin{equation}
dX(t)=\xi(t,X(t))dt+\sigma (t,X(t))dW(t),\quad 0<t<1,
\end{equation}
where $W(t)$ denotes a $d$-dimensional Brownian motion defined on a probability space (see Theorem \ref{thm7} in section \ref{sec:4}).
Let ${\cal P} (\mathbb{R}^d )$ denote the set of all Borel probability measures on $\mathbb{R}^d$.
For any $P_0, P_1\in {\cal P} (\mathbb{R}^d )$, there exists a unique product measure $\nu_0(dx)\nu_1(dy)$
that satisfies the following:
\begin{eqnarray}
P_1(dy)&=&\nu_1(dy)\int_{\mathbb{R}^d}p(0,x,1,y)\nu_0(dx),\label{r1.3}\\
P_0(dx)&=&\nu_0(dx)\int_{\mathbb{R}^d}p(0,x,1,y)\nu_1(dy).\label{r1.2}
\end{eqnarray}
This is Euler's equation of Schr\"odinger's problem and 
is called Schr\"odinger's functional equation or the Schr\"odinger system
(see \cite{1-S1,1-S2} and also \cite{1-Jamison1974} and Proposition 2.1 in \cite{2019OJM}).
Under some assumptions on $\sigma$ and $\xi$ (see, e.g.  (A.5)-(A.6) in section \ref{sec:4}), 
if $P_1(dy)\ll dy$, then there exists a unique weak solution $\{Y(t)\}_{0\le t\le 1}$ to the following
(see \cite{1-Jamison1975}) :
\begin{eqnarray}
dY(t)&=&\{a(t,Y(t))D_y\log h(t,Y(t))+\xi(t,Y(t))\}dt\label{r1.4}\\
&&\qquad +\sigma (t,Y(t))dW(t),\quad 0<t<1,\notag\\ 
P^{Y(0)}&=&P_0,\label{r1.5}
\end{eqnarray}
where $a(t,x):=\sigma(t,x)\sigma(t,x)^*$, $\sigma(t,x)^*$ denotes the transpose of $\sigma(t,x)$, 
$D_y:=\left(\partial/\partial y_i\right)_{i=1}^d$, 
$$h(t,x):=\int_{\mathbb{R}^d}p(t,x,1,y)\nu_1(dy),\quad 0\le t<1, x\in \mathbb{R}^d,$$
and $P^{Y(0)}$ denotes the probability law of $Y(0)$.
Besides, the following holds:
\begin{equation}
P^{(Y(0),Y(1))}(dxdy)=\nu_0(dx)p(0,x,1,y)\nu_1(dy),
\end{equation}
which implies that $P^{Y(1)}=P_1$ from (\ref{r1.3}).
$\{Y(t)\}_{0\le t\le 1}$ is called the h-path process for $\{X(t)\}_{0\le t\le 1}$ 
from two endpoint marginals  $P_0, P_1$ at $t=0,1$, respectively.
\begin{remark}
Schr\"odinger's functional equation (\ref{r1.3})-(\ref{r1.2}) is equivalent to the following:
\begin{eqnarray}
\overline h(1,y)&=&\int_{\mathbb{R}^d}p(0,x,1,y)
\left\{\int_{\mathbb{R}^d}p(0,x,1,z)\overline h(1,z)^{-1}P_1(dz)\right\}^{-1}
P_0(dx),\qquad\label{r1.6}
\\
\nu_ 1(dy)&:=&\overline h(1,y)^{-1}P_1(dy),\\
\nu_0(dx)&:=&\left\{\int_{\mathbb{R}^d}p(0,x,1,z)\nu_ 1(dz)\right\}^{-1}P_0(dx).
\end{eqnarray}
In particular, one only has to find a solution $\overline h(1,\cdot)$ in (\ref{r1.6}).
\end{remark}


Motivated by E. Schr\"odinger's quantum mechanics,
E. Nelson proposed the  problem of the construction of a Markov diffusion process from the Fokker-Planck equation.
We describe it.
Let $a$ and $b$ be, respectively, a $d\times d$ symmetric nonnegative definite matrix-valued 
and an $\mathbb{R}^d$-valued function on $[0,1]\times \mathbb{R}^d$,
and let $\{ P_t\}_{0\le t\le 1}\subset {\cal P} (\mathbb{R}^d )$. 
By $(a,b)\in {\bf A} (\{P_t\}_{0\le t\le 1})$, we mean that 
$a,b\in L^1([0,1]\times \mathbb{R}^d, dtP_t(dx))$ and
the following Fokker-Planck equation holds:
 for any $f \in C^{1,2}_b ([0,1]\times \mathbb{R}^d)$ and $t\in [0,1]$,
\begin{eqnarray}\label{1.1}
&&\int_{\mathbb{R}^d}f (t,x)P_t(dx)-\int_{\mathbb{R}^d}f(0,x)P_0(dx)\\
&=&\int_0^t ds\int_{\mathbb{R}^d}\biggl(\partial_sf(s,x)+\frac{1}{2}
\langle a(s,x), D_x^2f(s,x)\rangle 
+ \langle b(s,x),D_x f (s,x)\rangle\biggr) P_s(dx).\notag
\end{eqnarray}
Here $\partial_s:=\partial /\partial s$, 
$D_x^2:=\left(\partial^2/\partial x_i\partial x_j\right)_{i,j=1}^d$,
$\langle x,y\rangle$ denotes the inner product of $x, y\in \mathbb{R}^d$ and 
\begin{equation*}
\langle A,B\rangle :=\sum_{i,j=1}^d A_{ij}B_{ij},\quad A=(A_{ij})_{i,j=1}^d, B=(B_{ij})_{i,j=1}^d
\in M(d,\mathbb{R}).
\end{equation*}
We also write $(a,b)\in {\bf A}_0 (\{P_t\}_{0\le t\le 1})$
if $a,b\in L^1_{loc}([0,1]\times \mathbb{R}^d, dtP_t(dx))$ and
(\ref{1.1}) holds for all $f \in C^{1,2}_0 ([0,1]\times \mathbb{R}^d)$.
\begin{remark}
For $\{ P_t\}_{0\le t\le 1}$ in (\ref{1.1}), 
${\bf A} (\{P_t\}_{0\le t\le 1})$ is not necessarily a singleton (see \cite{1-CL1}-\cite{1-CL4}, \cite{1-MikamiCMP}).

\end{remark}

The following is a generalized version of Nelson's problem (see \cite{1-Nagasawa,1-NelsonBrownian,1-NelsonQuantum}).

\begin{definition}[Nelson's problem\label{def1}]
For any $\{ P_t\}_{0\le t\le 1}\subset {\cal P} (\mathbb{R}^d )$ such that ${\bf A}_0 (\{P_t\}_{0\le t\le 1})$ is not empty and for any $(a,b)\in {\bf A}_0 (\{P_t\}_{0\le t\le 1})$, construct 
a $d\times d$ matrix-valued function $\sigma(t,x)$ on $[0,1]\times \mathbb{R}^d$ and 
a semimartingale $\{X(t)\}_{0\le t\le 1}$
such that the following holds: for $(t,x)\in [0,1]\times\mathbb{R}^d$,
\begin{eqnarray}
a(t,x)&=&\sigma(t,x)\sigma(t,x)^*,\quad dtP_t(dx)\hbox{-a.e.},\label{1.2}\\
X(t)&=&X(0)+\int_0^t b(s,X(s))ds+\int_0^t \sigma (s, X(s))dW_X(s),\label{1.3}\\ 
P^{X(t)}&=&P_t.\label{1.4}
\end{eqnarray}
Here 
$W_X$ 
denotes a $d$-dimensional  
Brownian motion.
\end{definition}

When $\sigma(t,x)$ is nondegenerate, $W_X$ can be taken to be an $({\cal F}_t^X)$-Brownian motion, where ${\cal F}_t^X$ denotes $\sigma [X(s):0\le s\le t]$.
Otherwise (\ref{1.3}) means that $X(t)-X(0)-\int_0^t b(s,X(s))ds$
 is a local martingale with a quadratic variational process $\int_0^t a (s, X(s))ds$ (see, e.g. \cite{2-IW}).

The first result on Nelson's problem was given by  E. Carlen when $a$ is an identity matrix  (see \cite{1-[1],1-12}, and also \cite{1-13,1-62,1-120} for different approaches).
We generalized it to the case with a variable diffusion matrix (see \cite{1-MikamiCMP}).
P. Cattiaux, C. L\'eonard extensively generalized it to the case where
the jump-type Markov processes are also considered (see  \cite{1-CL1}-\cite{1-CL4}).
In these papers, they assumed the following condition.

\begin{definition}[Finite energy condition (FEC)\label{def2}]

\noindent
There exists $(a,b) \in {\bf A} (\{P_t\}_{0\le t\le 1})$
such that the following holds:
\begin{equation}\label{1.5}
\int_0^1 dt\int_{\mathbb{R}^d}\langle a(t,x)^{-1}b(t,x),b(t,x)\rangle P_t(dx)<\infty.
\end{equation}
\end{definition}

We describe a class of  stochastic optimal transportation problems  (SOTPs for short)
and  approaches to the h-path process and Nelson's problem by the SOTPs.

Fix a Borel measurable $d\times d$-matrix function $\sigma(t,x)$.
Let ${\cal A}$ denote the set of all $\mathbb{R}^d$-valued, continuous semimartingales
$\{ X(t)\}_{0\le t\le 1}$ on a (possibly different) complete filtered probability space
such that there exists a Borel measurable $\beta_X :[0,1]\times C([0,1])\longrightarrow \mathbb{R}^d$
for which the following holds:

\medskip
\noindent 
(i) $\omega \mapsto \beta_X (t,\omega )$ is ${\bf B}(C([0,t]))_+$-measurable for all $t\in [0,1]$,

\medskip
\noindent 
(ii)  $X(t)=X(0)+\int_0^t \beta_X (s,X)ds+\int_0^t \sigma (s,X(s))dW_X(s)$, $0\le t\le 1$,

\medskip
\noindent 
(iii)  
\begin{equation*}
E\left[\int_0^1\left\{ |\beta_X (t,X)|+|\sigma (t,X(t))|^2\right\}dt\right]<\infty.
\end{equation*}
Here ${\bf B}(C([0,t]))$ and ${\bf B}(C([0,t]))_+$ denote the Borel $\sigma$-field of $C([0,t])$
and $\cap_{s> t}{\bf B}(C([0,s]))$, respectively (see, e.g. \cite{2-LS}).
$|\cdot |:=\langle \cdot, \cdot\rangle^{1/2}$.

Let $L:[0,1]\times \mathbb{R}^d \times \mathbb{R}^d\longrightarrow [0,\infty )$ be continuous.
The following is a class of  the SOTPs
(see \cite{1-21-1,1-MTdual}, and also \cite{1-MikamiCMP,1-Mikami2002ECP,Mikami}).

\begin{definition}[Stochastic optimal transportation problems\label{def3}]

\noindent
 (1) For $P_0$, $P_1\in {\cal P} (\mathbb{R}^d )$,
\begin{equation}\label{1.6}
V (P_0,P_{1}):=\inf _{{\scriptstyle X\in {\cal A},} 
\atop{{\scriptstyle P^{X (t)}=P_t, t=0,1}}}\
E\biggl[\int_0^1 L(t,X (t);\beta_X (t,X))dt \biggr].
\end{equation}

\noindent
(2) For $\{ P_t\}_{0\le t\le 1}\subset {\cal P} (\mathbb{R}^d )$,
\begin{equation}\label{1.7}
{\bf V} (\{ P_t\}_{0\le t\le 1} ):=\inf _{{\scriptstyle X\in {\cal A},} 
\atop{{\scriptstyle P^{X (t)}=P_t, 0\le t\le 1}}}
E\biggl[\int_0^1 L(t,X (t);\beta_X (t,X))dt \biggr].
\end{equation}
If the set over which the infimum is taken is empty, then we set the infimum for infinity.
\end{definition}

Suppose that one knows the marginal probability distributions of a stochastic system at times $t=0, 1$ or $t\in [0,1]$.
To study the stochastic system on $[0,1]$ from the viewpoint of  the  principle of least action, one has to consider these kinds of problems.

\begin{remark}\label{rk2}
(i) The sets of stochastic processes over which the infimum are taken in (\ref{1.6})-(\ref{1.7}) can  be empty.
If $P_1(dx)\ll dx$, then the case when it is not empty is known
for  (\ref{1.6}) in \cite{1-Jamison1975} and for (\ref{1.7}) in \cite{1-Bog}, \cite{1-[1],1-12}, \cite{1-13}-\cite{1-CL4}, \cite{1-MikamiCMP}, \cite{1-MikamiMC}-\cite{1-Mikami2002ECP}, \cite{1-MikamiSemi,1-21-1,1-62,1-Tre,1-120}.
(ii) For $\{X(t)\}_{0\le t\le 1}\in {\cal A}$,
\begin{eqnarray}\label{1.8}
&&(\sigma(t,x)\sigma(t,x)^*, b_X(t,x):=E[\beta_X(t,X)|(t, X(t)=x)])_{(t,x)\in [0,1]\times \mathbb{R}^d}\qquad\\
&\in &{\bf A} (\{P^{X(t)}\}_{0\le t\le 1}).\notag
\end{eqnarray}
Indeed, by It\^o's formula, (\ref{1.1}) with $a=\sigma\sigma^*, b=b_X$ holds and by Jensen's inequality,
\begin{equation}\label{1.9J}
E[|b_X(t,X(t))|]=E[|E[\beta_X(t,X)|(t, X(t))]|]\le E[|\beta_X(t,X)|].
\end{equation}
\end{remark}

Schr\"odinger's problem which is a typical example of the SOTP
is $V_S:=V$ in (\ref{1.6}) when 
the following holds:
\begin{equation}\label{1.12}
L=\frac{1}{2}|\sigma(t,x)^{-1} (u-\xi (t,x))|^2
\end{equation}
(see, e.g. \cite{1-Leo2,Mikami,1-RT}).
If $V_S (P_0, P_1)$ is finite for $P_0$, $P_1\in {\cal P} (\mathbb{R}^d )$
and if $\sigma$ and $\xi$ satisfy nice conditions,
then the minimizer uniquely exists and is the h-path process with two endpoint marginals $P_0, P_1$
 in (\ref{r1.4})-(\ref{r1.5}) (see \cite{1-Dai,1-Follmer,Mikami,1-MTdual,1-RR,1-Zambrini}).

By  the continuum limit of $V (\cdot,\cdot)$, we considered Nelson's Problem 
in a more general setting, including the following case (see \cite{1-MikamiCMP,1-MikamiSemi}).

\begin{definition}[Generalized finite energy condition (GFEC)\label{def4}]

\noindent
There exists  $\gamma >1$ and $(a,b) \in {\bf A} (\{P_t\}_{0\le t\le 1})$ such that  
the following holds:
\begin{equation}\label{1.9}
\int_0^1 dt\int_{\mathbb{R}^d}\langle a(t,x)^{-1}b(t,x),b(t,x)\rangle ^{\frac{\gamma}{2}} P_t(dx)<\infty.
\end{equation}
\end{definition}

As an application of the Duality Theorem for ${\bf V}$, we also gave an approach to Nelson's Problem under the condition which includes the  GFEC 
 (see \cite{1-21-1}).

If (\ref{1.2})-(\ref{1.4}) hold, then they also say that the superposition principle holds.
When $a\equiv 0$, the superposition principle was studied in \cite{1-Amb,1-AmbTre,1-Mikami2000,1-Mikami2002ECP}.
D. Trevisan's result \cite{1-Tre} almost completely solved Nelson's problem (see also \cite{1-Bog,1-Figalli}).
In the case where the linear operator with the second order differential operator and  with the L\'evy measure is considered, it was studied in \cite{1-CL4,1-Rockner}.


\begin{theorem}[see \cite{1-Tre}\label{thm2.3}]
Suppose that there exists $\{ P_t\}_{0\le t\le 1}\subset {\cal P} (\mathbb{R}^d )$ such that 
$(a,b)\in {\bf A} (\{P_t\}_{0\le t\le 1})$ exists.
Then Nelson's problem (\ref{1.2})-(\ref{1.4}) has a solution.
\end{theorem}

In his problem, E. Nelson considered the case where $a=Identity$ and $b=D_x \psi (t,x)$ for some function $\psi$.
It turned out that the Nelson process is the minimizer of ${\bf V}_N:={\bf V}$ when (\ref{1.12}) with 
$\sigma=Identity$ and $\xi=0$ and the FEC hold (see Proposition 3.1 in \cite{1-MikamiCMP} and also Theorem \ref{thm1-4} in section 2).
Indeed, if $(a, D_x \psi_i)\in {\bf A} (\{P_t\}_{0\le t\le 1})$, $i=1,2$, then $D_x \psi_1=D_x \psi_2$,
$dtP_t(dx)$-a.e..
In this sense, we consider that Nelson's problem is the studies of the superposition principle and 
of the minimizer of ${\bf V}$.
In particular, if the superposition principle holds, then the set over which the infimum is taken 
in ${\bf V}$ is not empty and then one can consider a minimizer of ${\bf V}$,
provided it is finite.
There was a different approach by showing Proposition \ref{pp1} in section 2 via the Duality Theorems in Theorems \ref{thm1} and \ref{thm1-4} in section 2
(see \cite{1-21-1} and also \cite{1-MikamiCMP,1-MikamiSemi}).
It is also generalized by the superposition principle and our previous approach to the first part of 
Nelson's problem is not useful anymore.

In section 2, we improve our previous results on the SOTPs with a convex cost function by the superposition principle in Theorem \ref{thm2.3}.
More precisely, we prove that the SOTPs are equivalent to variational problems 
for probability measures given by the Fokker-Planck equation and to
those by a relaxed version of the Fokker-Planck equation (see Proposition \ref{pp1} in section 2).
In particular, we can prove the convexity and the lower-semicontinuity of the  SOTPs in marginal distributions
by a finite-dimensional approach though  the SOTPs are  variational problems for semimartingales.
It gives a new insight into the SOTPs and lets us revisit them.


In section 3, in the case where $d=1$ and where $a$ is not fixed,
we consider  slightly relaxed versions of the SOTPs of which cost functions are not supposed to be convex.
In this case, we need a generalization of D. Trevisan's result which was recently obtained by V. I. Bogachev, M. R\"ockner, S. V. Shaposhnikov.

\begin{theorem}[see \cite{1-Bog}\label{thm5}]
Suppose that there exists $\{ P_t\}_{0\le t\le 1}\subset {\cal P} (\mathbb{R}^d )$ such that 
$(a,b)\in {\bf A}_0 (\{P_t\}_{0\le t\le 1})$ exists and that the following holds:
\begin{equation}\label{3.11}
\int_0^1 dt\int_{\mathbb{R}^d}\frac{|a(t,x)|+|\langle x,b(t,x)\rangle|}{1+|x|^2}P_t(dx)<\infty.
\end{equation}
Then Nelson's problem (\ref{1.2})-(\ref{1.4}) has a solution.
\end{theorem}

As a fundamental problem of the stochastic optimal control theory,  
the test of the Markov property of a minimizer is known.
We also discuss this problem for a finite-time horizon stochastic optimal control problem.

In section 4, we study the semiconcavity and the  Lipschitz continuity of Schr\"odinger's problem $V_S$.

\section{SOTPs with a convex cost}\label{sec:2}
In this section, we discuss applications of D. Trevisan's result to
the Duality Theorems for the SOTPs in the case where $u\mapsto L(t,x;u)$ is convex and where $\sigma$ and $a=\sigma\sigma^*$ in (\ref{1.2}) are fixed.
We write $b\in {\bf A} (\{P_t\}_{0\le t\le 1})$ if $(a,b)\in {\bf A} (\{P_t\}_{0\le t\le 1})$ for the sake of simplicity  (see (\ref{1.1}) for notation).

As a preparation, we introduce two classes of marginal problems which play crucial roles
in the proof of the Duality Theorems for the SOTPs (see \cite{1-MikamiSemi,1-21-1}) and
which will be proved to be equivalent to the SOTPs by D. Trevisan's result.

The following can be considered as versions of the SOTPs for a flow of marginals which satisfy (\ref{1.1}).

\begin{definition}[SOTPs for marginal flows\label{def5}]

\noindent
(1) For $P_0$, $P_{1}\in {\cal P}(\mathbb{R}^d)$,
\begin{equation}\label{2.1}
{\rm v}(P_0,P_{1})
:=\inf _{{\scriptstyle b\in {\bf A} (\{Q_t\}_{0\le t\le 1}),}
\atop{{\scriptstyle Q_t=P_t, t=0,1}}}
\int_0^1 dt \int_{{\bf R }^d}L(t,x;b(t,x))Q_t(dx).
\end{equation}
(2) For $\{ P_t\}_{0\le t\le 1}\subset {\cal P}(\mathbb{R}^d)$,
\begin{equation}\label{2.2}
{\bf v} (\{ P_t\}_{0\le t\le 1} )
:=\inf_{b\in {\bf A} (\{ P_t\}_{0\le t\le 1})} \int_0^1 dt\int_{\mathbb{R}^d} L(t,x;b(t,x))P_t(dx).
\end{equation}
\end{definition}

For $\mu (dxdu)\in {\cal P} ( \mathbb{R}^d\times \mathbb{R}^d )$,
\begin{equation}\label{2.3}
\mu_{1}(dx):=\mu (dx\times \mathbb{R}^d),\quad \mu_{2}(du):=\mu (\mathbb{R}^d\times du).
\end{equation}
We write $\nu (dtdxdu)\in \tilde{\cal A}$ if the following holds.
(i) $\nu\in {\cal P} ([0,1]\times \mathbb{R}^d \times \mathbb{R}^d )$ and 
\begin{equation}
\int_{[0,1]\times \mathbb{R}^d\times \mathbb{R}^d} (|a(t,x)|+|u|)\nu (dtdxdu)<\infty.
\end{equation}
(ii)
$\nu (dtdxdu)=dt\nu (t,dxdu)$,
$\nu (t,dxdu)\in {\cal P} ( \mathbb{R}^d\times \mathbb{R}^d )$, 
$\nu_{1}(t, dx), \nu_{2}(t, du)\in {\cal P} ( \mathbb{R}^d)$, $dt-$a.e. and
$t\mapsto \nu_{1}(t, dx)$ has a weakly continuous version $\nu_{1,t}(dx)\in {\cal P} ( \mathbb{R}^d)$  for which
the following holds:
for any $t\in [0,1]$ and $f\in C^{1,2}_b ([0,1]\times \mathbb{R}^d)$,
\begin{eqnarray}\label{2.4}
&&\int_{\mathbb{R}^d}f (t,x)\nu_{1,t}(dx)-\int_{\mathbb{R}^d}f(0,x)\nu_{1,0}(dx)\\
&=&\int_{[0,t]\times \mathbb{R}^d\times \mathbb{R}^d} {\cal L}_{s,x,u} f(s,x)\nu (dsdxdu).\notag
\end{eqnarray}
Here  
\begin{equation}\label{2.5}
{\cal L}_{s,x,u} f(s,x):=\partial_s f(s,x)+\frac{1}{2}\langle a(s,x),D_x^2 f(s,x)\rangle +\langle u,D_x f (s,x)\rangle.
\end{equation}
We introduce a relaxed version of the problem above (see \cite{2-Gomes} and references therein
for related topics).

\begin{definition}[SOTPs for marginal measures\label{def6}]

\noindent
(1) For $P_0$, $P_{1}\in {\cal P}(\mathbb{R}^d)$,
\begin{equation}\label{2.6}
\tilde {\rm v}(P_0,P_{1})
:=\inf _{{\scriptstyle \nu \in \tilde{\cal A},}
\atop{{\scriptstyle \nu_{1,t}=P_t, t=0,1}}}
\int_{[0,1]\times \mathbb{R}^d\times \mathbb{R}^d}L(t,x;u)\nu (dtdxdu).
\end{equation}
(2) For $\{ P_t\}_{0\le t\le 1}\subset {\cal P}(\mathbb{R}^d)$,
\begin{equation}\label{2.7}
\tilde {\bf v} (\{ P_t\}_{0\le t\le 1} )
:=\inf _{{\scriptstyle \nu\in \tilde{\cal A},}
\atop{{\scriptstyle \nu_{1,t}=P_t, 0\le t\le 1}}}
\int_{[0,1]\times \mathbb{R}^d\times \mathbb{R}^d}L(t,x;u)\nu (dtdxdu).
\end{equation}
\end{definition}

\begin{remark}\label{rk3}
If $b\in {\bf A}(\{ P_t\}_{0\le t\le 1})$ and $X\in {\cal A}$, then 
$dtP_t(dx)\delta_{b(t,x)}(du)\in \tilde{\cal A}$
and  $dtP^{(X(t),\beta_X(t,X))}(dxdu)\in \tilde{\cal A}$, respectively.
Here $\delta_x$ denotes the delta measure on $\{x\}$.
In particular, $dtP^{(X(t),\beta_X(t,X))}(dxdu)$ is the distribution of a $[0,1]\times \mathbb{R}^d\times \mathbb{R}^d$-valued random variable $(t,X(t),\beta_X(t,X))$.
This is why we call (\ref{2.6})-(\ref{2.7}) SOTPs for marginal measures
(see also Lemma \ref{lm1} given later).
One can also identify $\{ P_t\}_{0\le t\le 1}\subset {\cal P}(\mathbb{R}^d)$ with 
$dtP_t(dx)\in  {\cal P}([0,1]\times\mathbb{R}^d)$ when ${\bf V}, {\bf v}$ and $\tilde {\bf v}$
are considered (see Theorem \ref{thm1-4} and also \cite{1-21-1,Mikami}).
\end{remark}

We introduce assumptions.

\noindent
(A.0.0). 
(i) $\sigma_{ij}\in C_b([0,1]\times \mathbb{R}^d)$, $i,j=1,\cdots ,d$.
(ii) $\sigma (\cdot)=(\sigma_{ij}(\cdot))_{i,j=1}^d$
is a nondegenerate $d\times d$-matrix function on $[0,1]\times \mathbb{R}^d$.

\noindent
(A.1). 
(i) $L\in C([0,1]\times \mathbb{R}^d \times \mathbb{R}^d;[0,\infty ))$.
(ii) $\mathbb{R}^d\ni u\mapsto L(t,x;u)$ is convex for $(t,x)\in [0,1]\times\mathbb{R}^d$.

\noindent
(A.2). 

\begin{equation*}
\lim_{|u|\to\infty} \frac{\inf\{L(t,x;u)|(t,x)\in [0,1]\times \mathbb{R}^d\}}
{|u|}=\infty.
\end{equation*}

The following proposition gives the relations among and the properties of three classes of 
the SOTPs stated in Definitions \ref{def3}, \ref{def5}, and \ref{def6} above.
In particular, it implies that they are equivalent in our setting and why they are all called the SOTPs.
It also implies the convexities and the lower semicontinuities of $V(P_0,P_{1})$ and $ {\bf V} (\{ P_t\}_{0\le t\le 1})$.

\begin{proposition}\label{pp1}
(i) Suppose that (A.1) holds.
Then the following holds:
\begin{equation}\label{2.8}
V(P_0,P_{1})={\rm v}(P_0,P_{1})=\tilde {\rm v}(P_0,P_{1}),\quad P_0, P_{1}\in {\cal P}(\mathbb{R}^d),
\end{equation}
\begin{equation}\label{2.9}
{\bf V} (\{ P_t\}_{0\le t\le 1} )= {\bf v} (\{ P_t\}_{0\le t\le 1} )=\tilde {\bf v} (\{ P_t\}_{0\le t\le 1} ), \{ P_t\}_{0\le t\le 1}\subset {\cal P}(\mathbb{R}^d).   
\end{equation}
(ii) Suppose, in addition, that (A.0.0,i) and (A.2) hold.
Then there exist minimizers $X$ of $V (P_0,P_{1})$ and $Y$ of ${\bf V} (\{ P_t\}_{0\le t\le 1} )$ for which 
\begin{equation}\label{2.10}
\beta_X (t,X)=b_X(t,X(t)),\quad \beta_Y (t,Y)=b_Y(t,Y(t)),
\end{equation}
provided $V (P_0,P_{1})$ and  ${\bf V} (\{ P_t\}_{0\le t\le 1} )$ are finite, respectively
(see (\ref{1.8}) for notation).

\noindent
(iii) Suppose, in addition, that (A.0.0,ii) holds and that $\mathbb{R}^d\ni u\mapsto L(t,x;u)$ is strictly convex for $(t,x)\in [0,1]\times\mathbb{R}^d$.
Then for any minimizers $X$ of $V (P_0,P_{1})$ and $Y$ of ${\bf V} (\{ P_t\}_{0\le t\le 1} )$,
(\ref{2.10}) holds and 
$b_X$ and $b_Y$ in (\ref{2.10}) are unique on the support of $dtP^{X(t)}(dx)$ and $dtP^{Y(t)}(dx)$,
respectively.
\end{proposition}

\begin{remark}\label{remark4}
Let $c\in C(\mathbb{R}^d\times \mathbb{R}^d;[0,\infty))$.
For $P_0, P_{1}\in {\cal P}(\mathbb{R}^d)$,
\begin{eqnarray}
&&T_M(P_0, P_{1}):=\inf\left\{\int_{\mathbb{R}^d}c(x,\varphi (x))P_0(dx)\biggl|P_0\varphi^{-1}=P_1\right\}\\
&\ge&T(P_0, P_{1}):=
\inf\left\{\int_{\mathbb{R}^d\times \mathbb{R}^d}c(x,y)\mu(dxdy)\biggl|
\mu_i=P_{i-1}, i=1,2\right\}\nonumber
\end{eqnarray}
(see (\ref{2.3}) for notation).
$T_M(P_0, P_{1})$ and $T(P_0, P_{1})$ are called Monge's and Monge-Kantorovich's problems, respectively.
The second equalities in (\ref{2.8})-(\ref{2.9}) are similar to the relation between Monge's and Monge-Kantorovich's problems 
since $\tilde {\rm v}$ and $\tilde {\bf v}$ are the infimums of linear functionals of measure 
(see, e.g. \cite{1-RR,1-Villani1}).
\end{remark}

Before we prove Proposition \ref{pp1}, we state its application to the SOTPs.

For any $s\ge 0$ and $P\in {\cal P} (\mathbb{R}^d )$, 
\begin{equation}\label{511}
{\bf \Psi}_{P}(s):= \biggl\{ \nu \in \tilde{\cal A}\biggl|\nu_{1,0}=P,
\int_{[0,1]\times \mathbb{R}^d\times\mathbb{R}^d}L(t,x;u)\nu (dtdxdu)\le s \biggl\}.
\end{equation}
Let ${\cal P} (\mathbb{R}^d )$ be endowed with a weak topology.
Then the following is known.

\begin{lemma}[see \cite{1-21-1}\label{lm1}]
Suppose that (A.0.0,i) and (A.1)-(A.2) hold. 
Then for any $s\ge 0$ and compact set $K\subset {\cal P} (\mathbb{R}^d )$, 
the set $\cup_{P\in K}{\bf \Psi}_{P}(s)$
is compact in $ {\cal P} ([0,1]\times \mathbb{R}^d \times \mathbb{R}^d )$.
\end{lemma}

Lemma \ref{lm1} was given in \cite{1-21-1} to prove the Duality Theorems for 
${\rm v}(P_0,P_{1})$ and ${\bf v}(\{ P_t\}_{0\le t\le 1})$.
By Proposition \ref{pp1}, it can be also used in the proof of 
the lower semicontinuities of ${V}(P_0,P_{1})$ and ${\bf V}(\{ P_t\}_{0\le t\le 1})$.
Besides, we do not need the following assumption anymore.

\noindent
(A).
\begin{equation}\label{2.12}
\Delta L(\varepsilon_1,\varepsilon_2):=\sup \frac{L(t,x;u)-L(s,y;u)}{1+L(s,y;u)} \to
0\quad\hbox{as $\varepsilon_1,$ $\varepsilon_2\downarrow 0$},
\end{equation}
where the supremum is taken over all $(t,x)$ and  $(s,y) \in [0,1]\times \mathbb{R}^d $
for which $|t-s|<\varepsilon_1$, $|x-y|<\varepsilon_2$ and over all $u\in \mathbb{R}^d $.

\noindent
This assumption can be used to prove the lower semicontinuity of the following (see \cite{2-IT}, Chapter 9.1):
\begin{equation}
AC([0,1];\mathbb{R}^d)\ni f\mapsto 
\int_0^1 L\left(t,f(t);\frac{df(t)}{dt}\right)dt.
\end{equation}

We state additional assumptions and the improved versions of the Duality Theorems for ${V}(P_0,P_{1})$ and ${\bf V}(\{ P_t\}_{0\le t\le 1})$.

\noindent
(A.0). 
$\sigma_{ij}\in C^{1}_b ([0,1]\times \mathbb{R}^d)$, $i,j=1,\cdots ,d$.

\noindent
(A.3). 
(i) $\partial_t L(t,x;u)$ and $D_x L(t,x;u)$ are bounded on 
$[0,1]\times \mathbb{R}^d \times B_R$ for all $R>0$, 
where $B_R:=\{ x\in \mathbb{R}^d ||x|\le R\}$.
(ii) $C_L$ is finite, where
\begin{equation}C_L:=\sup \left\{
\frac{L(t,x;u)}{1+L(t,y;u)}\biggr|0\le t\le 1, x, y, u \in\mathbb{R}^d\right\}.
\end{equation}

\begin{equation}\label{2.14}
H(t,x;z):=\sup\{\langle z,u\rangle -L(t,x;u)|u\in \mathbb{R}^d\}.
\end{equation}
The following is a generalization of \cite{1-21-1}, 
 in that we do not need the nondegeneracy of $a$ and the assumption (A)
and can be proved almost in the same way as in \cite{1-21-1}
 by  Proposition \ref{pp1} and by Lemma \ref{lm1}.
Indeed, in our previous papers, by the nondegeneracy of $a$, we made use of the Cameron-Martin-Maruyama-Girsanov formula to prove the convexity of $P\mapsto V(P_0,P)$, which we can avoid by Proposition \ref{pp1}.
The lower semicontinuity of $P\mapsto V(P_0,P)$
can be proved by  Proposition \ref{pp1} and by Lemma \ref{lm1}.
In \cite{2-TT}, they considered a similar problem and used a general property on the convex combination of probability measures on an enlarged space, which allows them not to assume the nondegeneracy of $a$,
though they assumed a condition which is similar to (A).

One can also find details in \cite{Mikami} (see \cite{GMT} for related topics).
We refer readers to \cite{2-CIL,1-FS,2-koike} on the viscosity solution.

\begin{theorem}[Duality Theorem for $V$\label{thm1}]
Suppose that (A.0)-(A.3) hold.
Then, for any $P_0$, $P_{1}\in {\cal P} (\mathbb{R}^d )$, 
\begin{eqnarray}\label{2.15}
&&V(P_0,P_1)={\rm v}(P_0,P_{1})=\tilde {\rm v}(P_0,P_{1})\\
&=&\sup_{f\in C_b^\infty (\mathbb{R}^d)}
\biggl\{\int_{\mathbb{R}^d}f(x)P_1 
(dx)-\int_{\mathbb{R}^d}\varphi (0,x;f)P_0(dx)\biggr\},\nonumber
\end{eqnarray}
where $\varphi (t,x;f)$ denotes the minimal bounded continuous viscosity solution 
to  the following HJB Eqn:
 on $[0,1)\times \mathbb{R}^d$,
\begin{eqnarray}\label{2.16}
\partial_t\varphi (t,x)+
\frac{1}{2}\langle a(t,x), D_x^2 \varphi (t,x)\rangle+H(t,x;D_x \varphi (t,x))&=&0,\\
\varphi (1,\cdot )&=&f.\notag
\end{eqnarray}
\end{theorem}




We introduce the following condition to replace $\varphi$ in (\ref{2.15}) by classical solutions to the HJB Eqn. (\ref{2.16}).

\noindent
(A.4). 
(i) ``$\sigma$ is an identity", or ``
$\sigma (\cdot)=(\sigma_{ij}(\cdot))_{i,j=1}^d$ is uniformly nondegenerate,
$\sigma_{ij}\in C^{1,2}_b ([0,1]\times \mathbb{R}^d)$, $i,j=1,\cdots ,d$, and
there exist functions $L_1$ and $L_2$ so that $L=L_1 (t,x)+L_2 (t,u)$".
(ii) $L(t,x;u)\in C^1([0,1]\times \mathbb{R}^d \times \mathbb{R}^d;[0,\infty ))$ and is strictly convex in $u$.
(iii) $L\in C^{1,2,0}_b ([0,1]\times \mathbb{R}^d\times B_R )$ for any $R>0$.

Since (A.4,i), (A.4,ii), and (A.4,iii) imply (A.0), (A.1), and (A.3,i), respectively,
the following holds from Theorem \ref{thm1}, in the same way as in \cite{1-21-1} (see also \cite{Mikami}).

\begin{corollary}\label{co1}
Suppose that 
(A.2), (A.3,ii), and (A.4) hold.
Then (\ref{2.15}) holds even if  the supremum is taken over all classical solutions
$\varphi\in C_b^{1,2} ([0,1]\times \mathbb{R}^d)$ to the HJB Eqn (\ref{2.16}).
Besides, for any $P_0, P_1\in  {\cal P}(\mathbb{R}^d)$ for which $V (P_0, P_1)$ is finite,
a minimizer $\{X(t)\}_{0\le t\le 1}$ of $V (P_0, P_1)$ exists and the following holds:
for any maximizing sequence $\{\varphi_n\}_{n\ge 1}$ of  (\ref{2.15}),
\begin{eqnarray}
0&=&\lim_{n\to\infty} E\biggl[\int_0^1 |L(t,X (t);\beta_X(t,X)) \\
&&\qquad-\{\langle\beta_X(t,X), D_x \varphi_n (t,X(t))\rangle -H(t,X (t);D_x\varphi_n(t,X(t)))\}|dt\biggr].\notag
\end{eqnarray}
In particular, there exists a subsequence $\{n_k \}_{k\ge 1}$ for which
\begin{equation}\label{2.38}
\beta_X(t,X) = \lim_{k\to\infty} D_z H(t,X(t);D_x \varphi_{n_k}(t,X(t))),\quad dtdP\hbox{{\rm -a.e.}}. 
\end{equation}
\end{corollary}

The following is also a generalization of \cite{1-21-1} 
and can be proved almost in the same way as in \cite{1-21-1} by  Proposition \ref{pp1} and Lemma \ref{lm1}. 

\begin{theorem}[Duality Theorem for ${\bf V}$\label{thm1-4}]
Suppose that (A.0)-(A.3) hold.
Then for any ${\bf P}:=\{ P_t\}_{0\le t\le 1}\subset{\cal P} (\mathbb{R}^d )$,
\begin{eqnarray}\label{227}
&&{\bf V} ({\bf P})={\bf v} ({\bf P})=\tilde {\bf v} ({\bf P} )\\
&=&\sup_{f\in C_b^\infty ([0,1]\times \mathbb{R}^d)}
\biggl\{\int_0^1\int_{\mathbb{R}^d}f(t,x)dtP_t (dx) 
-\int_{\mathbb{R}^d}\phi (0,x;f)P_0(dx)\biggr\},\nonumber
\end{eqnarray}
where 
$\phi (t,x;f)$ denotes the minimal bounded continuous viscosity solution of the following HJB Eqn:
on $[0,1)\times\mathbb{R}^d$,
\begin{eqnarray}\label{226}
\partial_t\phi (t,x)+\frac{1}{2}\langle a(t,x), D_x^2 \phi (t,x)\rangle+H(t,x;D_x \phi (t,x))+f(t,x)&=&0,
\qquad\\
\phi (1,x)&=&0.\notag
\end{eqnarray}
Suppose that (A.4) holds instead of (A.0), (A.1), and (A.3,i).
Then (\ref{227}) holds even if  the supremum is taken over all classical solutions
$\phi\in C^{1,2}_b ([0,1]\times \mathbb{R}^d)$ to the HJB Eqn (\ref{226}).
Besides,  if ${\bf V} ({\bf P})$ is finite, then a minimizer $\{X(t)\}_{0\le t\le 1}$ of ${\bf V} ({\bf P})$
exists and the following holds:
for any maximizing sequence $\{\phi_n\}_{n\ge 1}$ of  (\ref{227}),
\begin{eqnarray}
0&=&\lim_{n\to\infty} E\biggl[\int_0^1 |L(t,X (t);\beta_X(t,X)) \\
&&\qquad-\{\langle\beta_X(t,X), D_x \phi_n (t,X(t))\rangle -H(t,X (t);D_x\phi_n(t,X(t)))\}|dt\biggr].\notag
\end{eqnarray}
In particular, there exists a subsequence $\{n_k \}_{k\ge 1}$ for which
\begin{equation}
\beta_X(t,X) = \lim_{k\to\infty} D_z H(t,X(t);D_x \phi_{n_k}(t,X(t))),\quad dtdP\hbox{{\rm -a.e.}}. 
\end{equation}
\end{theorem}

\begin{remark}{\rm (see \cite{1-21-1,Mikami}\label{rk2-5})}
(i) Suppose that  (A.0)-(A.3) hold.
Then for any $f\in UC_b (\mathbb{R}^d)$,
the following is the minimal bounded continuous viscosity solution of 
the HJB equation (\ref{2.16}):
\begin{equation}\label{032.28}
\varphi (t,x;f)=
\sup_{X\in {\cal A}_t, X (t)=x}E\biggl[f(X(1))-\int_t^1 L(s, X(s);\beta_X(s,X))ds\biggl],
\end{equation}
where ${\cal A}_t$ denotes ${\cal A}$ with a time interval $[0,1]$ replaced by $[t,1]$.
(ii) Suppose that (A.0)-(A.3) with $L$ replaced by $L(t,x;u)-f(t,x)$ hold.
Then the following is the minimal bounded continuous viscosity solution of  the HJB equation (\ref{226}):
\begin{equation}\label{2.29}
\phi (t,x;f)=
\sup_{{\scriptstyle X\in {\cal A}_t,}\atop{{\scriptstyle X (t)=x}}}
E\biggl[\int_t^1 \left\{f(s,X(s))-L(s, X(s);\beta_X(s,X))\right\}ds\biggl].
\end{equation}
\end{remark}

We consider Schr\"odinger's and Nelson's problems, i.e.,  $V_S$ and $V_N$.
We introduce a new assumption.

\noindent
{\bf (A.4)'}. (\ref{1.12}) holds, $\sigma (\cdot)=(\sigma_{ij}(\cdot))_{i,j=1}^d$
is uniformly nondegenerate, and $a\in  C^{1,2}_b ([0, 1] \times \mathbb{R}^d;M(d,\mathbb{R})),
\xi \in  C^{1,2}_b ([0, 1] \times \mathbb{R}^d;\mathbb{R}^d)$.

(A.4)' implies 
(A.0)-(A.3).
Besides, for $f\in C^3_b (\mathbb{R}^d)$ and $f\in C^{1,2}_b ([0,1]\times \mathbb{R}^d)$,
the HJB equations (\ref{2.16}) and  (\ref{226}) have unique classical solutions  in $C^{1,2}_b ([0,1]\times \mathbb{R}^d)$, respectively.
They are also the minimal bounded continuous viscosity solutions of (\ref{2.16}) and  (\ref{226}),
respectively,
since they have the same representation formulas given in Remark \ref{rk2-5}
(see, e.g. \cite{2-Friedman,1-FS} on classical solutions and Lemma 4.5 in \cite{1-21-1} on viscosity solution).
In particular, the following holds though 
(A.4)' does not imply (A.4).

\begin{corollary}\label{co2}
Suppose that 
(A.4)' holds.
Then the assertions in Corollary \ref{co1}	 and Theorem \ref{thm1-4} hold.
\end{corollary}

\begin{remark}
If (\ref{1.12}) holds, then
\begin{equation}
L(t,x;u)- \{\langle u, z\rangle -H(t,x;z)\}
=\frac{1}{2}|\sigma (t,x)^{-1}(a(t,x)z-u+\xi(t,x))|^2.
\end{equation}
\end{remark}

In the rest of this section, we prove Proposition \ref{pp1}.

(Proof of Proposition \ref{pp1}).
We prove (i).
For $\{X(t)\}_{0\le t\le 1}\in {\cal A}$, by Jensen's inequality,
\begin{equation}\label{2.20}
E\biggl[\int_0^1 L(t,X (t);\beta_X (t,X))dt \biggr]
\ge E\biggl[\int_0^1 L(t,X (t);b_X(t,X(t)))dt \biggr].
\end{equation}
Theorem \ref{thm2.3} implies the first equalities of (\ref{2.8})-(\ref{2.9})  (see Remark \ref{rk2}, (ii)).
For $\nu\in \tilde{\cal A}$, 
\begin{equation}\label{2.21}
b_\nu (t,x):=\int_{\mathbb{R}^d}u\nu (t,x,du),
\end{equation}
where $\nu (t,x,du)$ denotes a regular conditional probability of $\nu$ given $(t,x)$.
Then by Jensen's inequality,
\begin{equation}\label{2.22}
\int_{[0,1]\times \mathbb{R}^d\times \mathbb{R}^d}L(t,x;u)\nu (dtdxdu)
\ge \int_0^1 dt\int_{\mathbb{R}^d}L(t,x;b_\nu (t,x))\nu_{1,t} (dx).
\end{equation}
$b_\nu \in {\bf A} (\{\nu_{1,t}\}_{0\le t\le 1})$ from (\ref{2.4}),
since by Jensen's inequality,
\begin{equation*}
\int_{[0,1]\times \mathbb{R}^d}|b_\nu (t,x)|dt\nu_{1,t}(dx)
\le \int_{[0,1]\times \mathbb{R}^d\times \mathbb{R}^d}|u|\nu(dtdxdu)<\infty,
\end{equation*}
and for any $t\in [0,1]$ and $f\in C^{1,2}_b ([0,1]\times \mathbb{R}^d)$,
\begin{eqnarray}\label{2.23}
&&\int_{[0,t]\times \mathbb{R}^d\times \mathbb{R}^d} \langle  u, D_x f(s,x)\rangle \nu (dsdxdu)\\
&=&\int_0^t ds \int_{\mathbb{R}^d} \langle b_\nu (s,x), D_x f(s,x)\rangle \nu_{1,s} (dx).\notag
\end{eqnarray}
This implies the second equalities of (\ref{2.8})-(\ref{2.9})  (see Remark \ref{rk3}).

\noindent
The proof  of (ii) is done by Lemma \ref{lm1}, (\ref{2.23}), and Theorem \ref{thm2.3}.

\noindent
We prove (iii).
From (\ref{2.20}) and the strict convexity of $u\mapsto L(t,x;u)$,  (\ref{2.10}) holds.
For $b\in {\bf A} (\{P_t\}_{0\le t\le 1})$, 
$P_t(dx)\ll dx$, $dt$-a.e. from (A.0.0,ii), since $a, b\in L^1([0,1]\times\mathbb{R}^d,dtP_t(dx))$
(see \cite{2-BKR}, p. 1042, Corollary 2.2.2).
For $\{p_{i}(t,x)dx\}_{0\le t\le 1}\subset{\cal P}(\mathbb{R}^d)$, $b_i\in {\bf A} (\{p_{i}(t,x)dx\}_{0\le t\le 1})$, $i=0,1$, and $\lambda\in [0,1]$,
\begin{equation}\label{2.27}
p_\lambda :=(1-\lambda )  p_0+\lambda p_1 ,\quad 
b_\lambda :=1_{(0,\infty)}(p_\lambda)\frac{(1-\lambda )p_0b_0+\lambda p_{1}b_1}{p_\lambda},
\end{equation}
where $1_A(x)$ denotes an indicator function of $A\subset \mathbb{R}$.
Then $b_\lambda\in {\bf A}(\{p_\lambda (t,x)dx\}_{0\le t\le 1})$ and
\begin{eqnarray}\label{2.28}
&&\int_0^1 dt \int_{{\bf R }^d}L(t,x;b_\lambda (t,x))p_\lambda (t,x)dx\\
&\le&(1-\lambda )\int_0^1 dt \int_{{\bf R }^d}L(t,x;b_0 (t,x))p_0(t,x)dx\notag\\
&&\qquad+\lambda \int_0^1 dt \int_{{\bf R }^d}L(t,x;b_1 (t,x))p_1 (t,x)dx.\notag
\end{eqnarray}
Here the equality holds if and only if $b_0=b_1$ $dtdx$-a.e. 
on the set $\{(t,x)\in [0,1]\times \mathbb{R}^d| p_0 (t,x)p_1 (t,x)>0\}$.
$\Box$

\section{Stochastic optimal transport with a nonconvex cost}

In this section, in the case where $d=1$  and where $a$ is not fixed,
we consider slightly relaxed versions of the SOTPs of which cost functions are not supposed to be convex.
As a fundamental problem of the stochastic optimal control theory,  
the test of the Markov property of a minimizer of a stochastic optimal control problem is known.
We also consider the Markov property of the minimizer of a finite-time horizon stochastic control problem.
Our previous result \cite{1-MikamiMC} proved it in a one-dimensional case by the optimal transportation problem with a concave cost.
We generalize it by Theorem \ref{thm5} in section 1.

Since $a$ is not fixed in this section, we consider a new class of semimartingales.
Let $u=\{u(t)\}_{0\le t\le 1}$ and $\{W(t)\}_{0\le t\le 1}$ 
be a progressively measurable real valued process and a one-dimensional Brownian motion 
on the same  complete filtered probability space, respectively.
The probability space under consideration is not fixed in this section.
Let $\sigma :[0,1]\times \mathbb{R}\longrightarrow \mathbb{R}$
be a Borel measurable function.
Let $Y^{u,\sigma}=\{Y^{u,\sigma}(t)\}_{0\le t\le 1}$ be  a continuous semimartingale such that
the following holds weakly:
\begin{equation}\label{3.1}
Y^{u,\sigma}(t)=Y^{u,\sigma}(0)+\int_0^t u (s)ds
+\int_0^t \sigma (s,Y^{u,\sigma}(s))dW(s), \quad 0\le t\le 1,
\end{equation}
provided it exists.

For $r> 0$,
\begin{equation}
{\cal U}_{r}:=\left\{ (u,\sigma) \biggl| 
E\left[\int_0^1 \left(\frac{\sigma(t,Y^{u,\sigma}(t))^2}{1+|Y^{u,\sigma}(t)|^2}+|u(t)|\right)dt\right]<\infty,|\sigma|\ge r\right\},
\end{equation}
\begin{equation}
{\cal U}_{r,Mar}:=\left\{ (u,\sigma) \in {\cal U}_{r}| 
u(\cdot)=b_{Y^{u,\sigma}}(\cdot,Y^{u,\sigma}(\cdot))\right\},
\end{equation}
where $b_{Y^{u,\sigma}}(t,Y^{u,\sigma}(t)):=E[u(t)|(t,Y^{u,\sigma}(t))]$.
For $(u,\sigma)\in {\cal U}_{r}$,
\begin{eqnarray}
F^{Y^{u, \sigma}}_t (x)&:=&P(Y^{u, \sigma}(t)\le x),\\ \label{3.2}
G^u_t (x)&:=&P(u(t)\le x),\\ \label{3.3}
\tilde b_{u, Y^{u, \sigma}}(t,x)&:=&(G^u_t)^{-1} (1-F^{Y^{u, \sigma}}_t (x)),\quad (t,x)\in [0,1]\times \mathbb{R}.\label{3.4}
\end{eqnarray}
Here for a distribution function $F$ on $\mathbb{R}$,
\begin{equation*}F^{-1}(v):=\inf\{x\in \mathbb{R}| F(x)\ge v\},\quad 0<v<1.
\end{equation*}
$F^{-1}$ is called the quasi-inverse of $F$ (see, e.g. \cite{1-64,1-RR,1-72}).
\begin{equation}\label{3.5}
p^{Y^{u, \sigma}}(t,x):=\frac{P^{Y^{u,\sigma}(t)}(dx)}{dx}
\end{equation}
exists $dt$-a.e.
since $r$ is positive and $(\sigma^2, b_{Y^{u,\sigma}})\in {\bf A}_0(\{P^{Y^{u,\sigma}(t)}\}_{0\le t\le 1})$
 (see \cite{2-BKR}, p. 1042, Corollary 2.2.2).
 Indeed, by Jensen's inequality,
\begin{equation}\label{3.10}
\int_{\mathbb{R}} |b_{Y^{u, \sigma}}(t,y)|p^{Y^{u, \sigma}}(t,y)dy
=E[|E[u(t)|(t,Y^{u, \sigma}(t))]|]\le E[|u(t)|].
\end{equation}
From the idea of covariance kernels (see \cite{3-8,3-9,3-52,3-56}),
\begin{eqnarray}\label{3.6}
&&\tilde a_{u,Y^{u, \sigma}}(t,x)\\
&:=&1_{(0,\infty )} (p^{Y^{u, \sigma}}(t,x))
\frac{2\int_{-\infty}^x (\tilde b_{u,Y^{u, \sigma}}(t,y)
-b_{Y^{u, \sigma}}(t,y))p^{Y^{u, \sigma}}(t,y)dy}{p^{Y^{u, \sigma}}(t,x)}.\nonumber
\end{eqnarray}
The following holds and will be proved later.

\begin{theorem}\label{thm3-5}
Let $r> 0$.
For $(u,\sigma)\in {\cal U}_{r}$,
there  exists $\tilde u$ such that
$(\tilde u,\tilde\sigma:=(\sigma^2+\tilde a_{u,Y^{u,\sigma}})^{\frac{1}{2}})\in {\cal U}_{r,Mar}$ and that the following holds:
\begin{eqnarray}
P^{Y^{\tilde u,\tilde\sigma}(t)}&=&P^{Y^{u,\sigma}(t)}, \quad t\in [0,1],\\
b_{Y^{\tilde u,\tilde\sigma}}&=&\tilde b_{u,Y^{u,\sigma}},\\
P^{b_{Y^{\tilde u,\tilde\sigma}}(t,Y^{\tilde u,\tilde\sigma}(t))}&=&P^{u(t)},\quad dt{\rm-a.e.}.
\end{eqnarray}
\end{theorem}

For $r> 0$  and $\{P_t\}_{0\le t\le 1}\subset{\cal P}(\mathbb{R})$,
\begin{eqnarray}
{\bf A}_{0,r}(\{P_t\}_{0\le t\le 1})
&:=&\biggl\{ (a,b) \in {\bf A}_{0}(\{P_t\}_{0\le t\le 1})\biggl| a\ge r^2,\\
&&\quad \int_0^1 dt\int_{\mathbb{R}^d}\left(\frac{a(t,x)}{1+|x|^2}+|b(t,x)|\biggl)P_t(dx)<\infty\right\}.\notag
\end{eqnarray}
Let $L_1, L_2:[0,1]\times \mathbb{R}\longrightarrow [0,\infty)$ be Borel measurable.
For $(u,\sigma)$,
\begin{equation}
J(u,\sigma):=E\biggl[\int_0^1 (L_1 (t,Y^{u,\sigma} (t))+L_2 (t,u(t)))dt\biggl].
\end{equation}
For $(a,b)\in {\bf A}_{0}(\{P_t\}_{0\le t\le 1})$,
\begin{equation}
I(\{P_t\}_{0\le t\le 1},a,b):=\int_0^1 dt\int_{\mathbb{R}^d}(L_1 (t,x)+L_2(t,b(t,x)))P_t(dx).
\end{equation}
One easily obtain the following from Theorems \ref{thm5} and \ref{thm3-5}.

\begin{corollary}\label{co3}
Suppose that $L_1, L_2:[0,1]\times \mathbb{R}\longrightarrow [0,\infty)$ are Borel measurable.
Then for any $r> 0$, the following holds.
(i) For any $P_0,P_1\in {\cal P}(\mathbb{R})$,
\begin{eqnarray}\label{3.12}
&&\inf\{J(u,\sigma)|
(u,\sigma)\in {\cal U}_{r}, P^{Y^{u,\sigma}(t)}=P_t, t=0,1\}\\
&=&\inf\{J(u,\sigma)|(u,\sigma)\in {\cal U}_{r, Mar}, P^{Y^{u,\sigma}(t)}=P_t, t=0,1\}
\notag\\
&=&\inf\{I(\{Q_t\}_{0\le t\le 1},a,b)| (a,b)\in {\bf A}_{0,r}(\{Q_t\}_{0\le t\le 1}),
Q_t=P_t, t=0,1\}.\notag
\end{eqnarray}
In particular,  if there exists a minimizer in (\ref{3.12}), 
then there exists a minimizer $(u,\sigma)\in {\cal U}_{r, Mar}$.
(ii) For any $\{P_t\}_{0\le t\le 1}\subset{\cal P}(\mathbb{R})$,
\begin{eqnarray}\label{3.121}
&&\inf\{J(u,\sigma)|
(u,\sigma)\in {\cal U}_{r}, P^{Y^{u,\sigma}(t)}=P_t, 0\le t\le 1\}\\
&=&\inf\{J(u,\sigma)|(u,\sigma)\in {\cal U}_{r, Mar}, P^{Y^{u,\sigma}(t)}=P_t, 0\le t\le 1\}
\notag\\
&=&\inf\{I(\{P_t\}_{0\le t\le 1},a,b)| (a,b)\in {\bf A}_{0,r}(\{P_t\}_{0\le t\le 1})\}.\notag
\end{eqnarray}
In particular,  if there exists a minimizer in (\ref{3.121}), 
then there exists a minimizer $(u,\sigma)\in {\cal U}_{r,Mar}$.
\end{corollary}

Suppose that $L:[0,1]\times\mathbb{R}\times \mathbb{R}\longrightarrow [0,\infty),
\Psi:\mathbb{R}\longrightarrow [0,\infty)$ are Borel measurable.
Then for any $P_0\in {\cal P}(\mathbb{R})$,
\begin{eqnarray}
&&\inf_{{\scriptstyle (u,\sigma)\in {\cal U}_{r},}
\atop{{\scriptstyle P^{Y^{u,\sigma}(0)}=P_0}}}
E\left[\int_0^1 L(t,Y^{u, \sigma}(t);u(t))dt+\Psi (Y^{u, \sigma} (1))\right]\\
&=&\inf_{P\in{\cal P}(\mathbb{R})}\left\{V_r(P_0,P)+\int_{\mathbb{R}} \Psi (x)
P(dx)\right\},\notag
\end{eqnarray}
where $V_r$ denotes $V$ with ${\cal A}$ replaced by $\{Y^{u, \sigma}|(u, \sigma)\in {\cal U}_{r}\}$. 
In particular, we easily obtain the following from Corollary \ref{co3}.

\begin{corollary}
In addition to the assumption of Corollary \ref{co3},
suppose that $\Psi:\mathbb{R}\longrightarrow [0,\infty)$ is Borel measurable.
Then for any $r> 0$ and $P_0\in {\cal P}(\mathbb{R})$,
\begin{eqnarray}\label{3.13}
&&\inf\{J(u,\sigma)+E[\Psi (Y^{u,\sigma}(1))]|
(u,\sigma)\in {\cal U}_{r}, P^{Y^{u,\sigma}(0)}=P_0\}\\
&=&\inf\{J(u,\sigma)+E[\Psi (Y^{u,\sigma}(1))]|(u,\sigma)\in {\cal U}_{r, Mar}, P^{Y^{u,\sigma}(0)}=P_0\}.\notag
\end{eqnarray}
In particular,  if there exists a minimizer in (\ref{3.13}), 
then there exists a minimizer $(u,\sigma)\in {\cal U}_{r, Mar}$.
\end{corollary}

We prove Theorem \ref{thm3-5} by Theorem \ref{thm5}.

(Proof of Theorem \ref{thm3-5}).
For $(u,\sigma)\in {\cal U}_{r}$, the following holds (see \cite{1-MikamiMC}):
\begin{equation}\label{624}
\tilde a_{u,Y^{u, \sigma}}(t,\cdot )\ge 0,
\quad P^{\tilde b_{u,Y^{u, \sigma}}(t,Y^{u, \sigma} (t))}=P^{u(t)},\quad dt{\rm -a.e.}.
\end{equation}
Indeed,
\begin{eqnarray*}
\int_{-\infty}^x \tilde b_{u,Y^{u, \sigma}}(t,y)p^{Y^{u, \sigma}}(t,y)dy
&=&E[(G^u_t)^{-1} (1-F^{Y^{u, \sigma}}_t (Y^{u, \sigma} (t)));Y^{u, \sigma} (t)\le x],\\
\int_{-\infty}^{x} b_{Y^{u, \sigma}}(t,y)p^{Y^{u, \sigma}}(t,y)dy
&=&E[E[u(t)|(t,Y^{u, \sigma}(t))];Y^{u, \sigma} (t)\le x]\\
&=&E[u(t);Y^{u, \sigma}(t)\le x].
\end{eqnarray*}
For an $\mathbb{R}^2$-valued random variable $Z=(X,Y)$ on a probability space, 
\begin{eqnarray*}
E[Y;X\le x]&=&\int_0^\infty \{F_X(x)-F_Z(x,y)\}dy-\int_{-\infty}^0 F_Z(x,y)dy,\\
F_Z(x,y)&\ge &\max (F_X(x)+F_Y(y)-1,0)\\
&=&P(F_X^{-1} (U)\le x,F_Y^{-1} (1-U)\le y),
\end{eqnarray*}
where $F_X$ denotes the distribution function of $X$ and $U$ is a uniformly distributed random variable on $[0,1]$.
 The distribution functions of $F_X^{-1} (U)$ and $F_Y^{-1} (1-U)$ are $F_X$ and $F_Y$, respectively.
From (\ref{3.5}), $F^{Y^{u, \sigma}}_t (Y^{u, \sigma}(t))$ is uniformly distributed on $[0,1]$ and $(F^{Y^{u, \sigma}}_t)^{-1} (F^{Y^{u, \sigma}}_t (Y^{u, \sigma} (t)))=Y^{u, \sigma}(t)$, $P$-a.s., $dt$-a.e. (see \cite{1-21} or, e.g. \cite{1-64,1-RR,1-72}).

It is easy to see that the following holds from (\ref{3.10}) and (\ref{624}): 
\begin{equation*}
(\sigma^2+\tilde a_{u,Y^{u,\sigma}}, \tilde b_{u,Y^{u,\sigma}})\in {\bf A}_0(\{P^{Y^{u,\sigma}(t)}\}_{0\le t\le 1}).
\end{equation*}
Indeed, from (\ref{624}), the following holds:
\begin{equation}\label{3.109}
E\left[\int_0^1|\tilde b_{u,Y^{u,\sigma}}(t,Y^{u,\sigma} (t))|dt\right]=E\left[\int_0^1|u(t)|dt\right]<\infty.
\end{equation}
The following will be proved below:
\begin{equation}\label{3.110}
E\left[\int_0^1 \frac{\tilde a_{u,Y^{u,\sigma}}(t,Y^{u,\sigma}(t))}{1+|Y^{u,\sigma}(t)|^2}dt\right]<\infty.
\end{equation}
(\ref{3.109})-(\ref{3.110}) and Theorem \ref{thm5} complete the proof.
We prove (\ref{3.110}).
\begin{eqnarray}
&&
E\left[\int_0^1 \frac{\tilde a_{u,Y^{u,\sigma}}(t,Y^{u,\sigma}(t))}{2(1+|Y^{u,\sigma}(t)|^2)}dt\right]\\
&=&\int_0^1dt\int_{\mathbb{R}}\frac{1}{1+x^2}dx\int_{-\infty}^x (\tilde b_{u,Y^{u,\sigma}}(t,y)
-b_{Y^{u,\sigma}}(t,y))p^{Y^{u,\sigma}}(t,y)dy.\notag
\end{eqnarray}
From (\ref{624}),
\begin{eqnarray}
\int_{-\infty}^\infty \tilde b_{u,Y^{u,\sigma}}(t,y)p^{Y^{u,\sigma}}(t,y)dy
&=&E[u(t)]=E[E[u(t)|(t,Y^{u,\sigma}(t))]]\\
&=&\int_{-\infty}^\infty b_{Y^{u,\sigma}}(t,y)p^{Y^{u,\sigma}}(t,y)dy,\quad dt\hbox{-a.e.}.\notag
\end{eqnarray}
In particular, the following holds $dt$-a.e.:
\begin{eqnarray}
&&\int_{\mathbb{R}}\frac{1}{1+x^2}dx\int_{-\infty}^x (\tilde b_{u,Y^{u,\sigma}}(t,y)
-b_{Y^{u,\sigma}}(t,y))p^{Y^{u,\sigma}}(t,y)dy\\
&=&\int_{-\infty}^0\frac{1}{1+x^2}dx\int_{-\infty}^x (\tilde b_{u,Y^{u,\sigma}}(t,y)
-b_{Y^{u,\sigma}}(t,y))p^{Y^{u,\sigma}}(t,y)dy\notag\\
&&-\int_0^{\infty}\frac{1}{1+x^2}dx\int_x^{\infty} (\tilde b_{u,Y^{u,\sigma}}(t,y)
-b_{Y^{u,\sigma}}(t,y))p^{Y^{u,\sigma}}(t,y)dy\notag\\
&=&\int_{-\infty}^0(\tilde b_{u,Y^{u,\sigma}}(t,y)
-b_{Y^{u,\sigma}}(t,y))p^{Y^{u,\sigma}}(t,y)dy\int_{y}^0\frac{1}{1+x^2}dx\notag\\
&&-\int_0^{\infty}(\tilde b_{u,Y^{u,\sigma}}(t,y)
-b_{Y^{u,\sigma}}(t,y))p^{Y^{u,\sigma}}(t,y)dy\int_0^y\frac{1}{1+x^2}dx\notag\\
&\le &\int_{\mathbb{R}}|\arctan y|( |\tilde b_{u,Y^{u,\sigma}}(t,y)|+|b_{Y^{u,\sigma}}(t,y)|)p^{Y^{u,\sigma}}(t,y)dy.\notag
\end{eqnarray}
Since $|\arctan y|$ is bounded,
(\ref{3.10}) and (\ref{3.109}) completes the proof of (\ref{3.110}).
$\Box$

\section{Semiconcavity and continuity of Schr\"odinger's Problem\label{sec:4}}

Proposition \ref{pp1} and Lemma \ref{lm1} imply that
${\cal P} (\mathbb{R}^d \times \mathbb{R}^d )\ni P\times Q\mapsto V(P,Q)$
 is convex and lower semicontinuous.
In this section, we give a sufficient condition under which 
for a fixed $Q\in {\cal P} (\mathbb{R}^d )$,
$L^2 (\Omega, P;\mathbb{R}^d)\ni X\mapsto V_S(P^X,Q)$
is semiconcave and is continuous (see (\ref{1.12}) for notation).
More precisely, we show that there exists $C>0$ such that for a fixed $Q\in {\cal P} (\mathbb{R}^d )$,
\begin{equation*}
L^2 (\Omega, P;\mathbb{R}^d)\ni X\mapsto V_S(P^X,Q)-CE[|X|^2]
\end{equation*}
is concave and is continuous.
Here  $L^2 (\Omega, P;\mathbb{R}^d)$ denotes
the space of all square integrable functions from  a probability space $(\Omega, {\cal F}, P)$ to  $(\mathbb{R}^d,{\bf B}(\mathbb{R}^d))$.
Let $W_2$ denote the Wasserstein distance of order 2, i.e. $T^{1/2}$ with $c=|y-x|^2$ in
Remark \ref{remark4}.
We also show the Lipschitz continuity of ${\cal P}_{2} (\mathbb{R}^d )\ni P\mapsto V_S(P,Q)$ 
in $W_2$ (see (\ref{4.5}) for notation).

We first describe the assumptions in this section.

\noindent
(A.5) $\sigma(t,x)=(\sigma^{ij}(t,x))_{i,j=1}^d$, $(t,x)\in [0,1]\times \mathbb{R}^d$, is a $d\times d$-matrix.
$a(t,x):=\sigma(t,x)\sigma(t,x)^*$, $(t,x)\in [0,1]\times \mathbb{R}^d$, is uniformly nondegenerate,
bounded, once continuously differentiable, and uniformly H\"older continuous.
$D_x a(t,x)$ is bounded and the first derivatives of $a(t,x)$ are uniformly H\"older continuous
in $x$ uniformly in $t\in [0,1]$.

\noindent
(A.6) $\xi(t,x):[0,1]\times \mathbb{R}^d\longrightarrow \mathbb{R}^d$ is bounded, continuous, and uniformly H\"older continuous
in $x$ uniformly in $t\in [0,1]$.
 
\begin{remark}
(A.5)-(A.6) imply (A.0.0), (A.1), and (A.2) for (\ref{1.12}).
(A.4)' implies (A.0)-(A.3) and (A.5)-(A.6).
\end{remark}

We describe the following fact.

\begin{theorem}\label{thm7}
Suppose that (A.5)-(A.6) hold.
Then for any $P_0\in {\cal P}(\mathbb{R}^d)$,
 the following SDE has the unique weak solution with a positive continuous transition probability density $p(t,x;s,y)$, $0\le t<s\le 1$, $x,y\in \mathbb{R}^d$:
\begin{eqnarray}\label{31}
d{\bf X}(t)&=&\xi(t,{\bf X}(t))dt+\sigma (t,{\bf X}(t))dW_{\bf X}(t),\quad 0< t<1,\\
P^{{\bf X} (0)}&=&P_0\notag 
\end{eqnarray}
(see \cite{1-Jamison1975}).
Besides, there exist constants $C_1, C_2>0$ such that 
\begin{equation}\label{4.11}
-C_1+C_2^{-1}|x-y|^2\le -\log  p(0,x;1,y)\le C_1+C_2|x-y|^2,\quad x,y\in \mathbb{R}^d
\end{equation}
 (see \cite{3-Aronson,2-Friedman}).
\end{theorem}

\begin{remark}
If $V_S(P,Q)$ is finite, 
 then the distribution of the minimizer $X$ of $V_S(P,Q)$ is absolutely continuous with respect to $P^{\bf X}$.
 In particular, $Q(dx)\ll dx$ under (A.5)-(A.6).  
 Indeed, $V_S(P,Q)$ is the relative entropy of $P^{X}$ with respect to $P^{\bf X}$
 and $P^{{\bf X}(1)}$ has a density
 (see the discussion below Remark \ref{rk2}).
\end{remark}

 We recall the definition of displacement convexity.
 
\begin{definition}[Displacement convexity (see \cite{3-McCann})]
Let $G:{\cal P}(\mathbb{R}^d )\longrightarrow \mathbb{R}\cup\{\infty\}$.
$G$ is displacement convex if the following is convex:
for any $\rho_0, \rho_1\in {\cal P}(\mathbb{R}^d )$ and 
convex function $\varphi:\mathbb{R}^d\longrightarrow \mathbb{R}\cup\{\infty\}$, 
\begin{equation}\label{3.35}
[0,1]\ni t\mapsto G(\rho_t),
\end{equation}
where $\rho_t:= \rho_0(id+t(D\varphi-id))^{-1}$, $0< t<1$,
provided $\rho_1= \rho_0(D\varphi)^{-1}$ and $\rho_t$ can be defined.
Here $id$ denotes an identity mapping.
\end{definition}

Recall that a convex function is differentiable $dx$-a.e. in the interior of its domain (see, e.g. \cite{1-Villani1}) and $\rho_t$ in (\ref{3.35}) is well defined if $\rho_0\in {\cal P}_{2,ac} (\mathbb{R}^d )$
and if $\rho_1\in {\cal P}_{2} (\mathbb{R}^d )$ (see, e.g. \cite{1-Villani1}).
Here
\begin{eqnarray}
{\cal P}_2 (\mathbb{R}^d )&:=&\left\{P\in {\cal P} (\mathbb{R}^d )
\biggl|\int_{\mathbb{R}^d}|x|^2P(dx)<\infty\right\},\label{4.5}\\
{\cal P}_{ac} (\mathbb{R}^d )&:=&\{p (x)dx\in {\cal P} (\mathbb{R}^d )\},\\
{\cal P}_{2,ac} (\mathbb{R}^d )&:=&{\cal P}_{2} (\mathbb{R}^d )\cap {\cal P}_{ac} (\mathbb{R}^d ).
\end{eqnarray}

The following implies that 
$L^2 (\Omega, P;\mathbb{R}^d)\ni X\mapsto  V_S(P^X,Q)$
 is semiconvave for a fixed $Q\in {\cal P}_{ac} (\mathbb{R}^d )$
 and will be proved later.

\begin{theorem}\label {pp34}
Suppose that  
(A.4)' holds and that there exists a constant $C>0$ such that
$x\mapsto \log p(0,x;1,y) +C|x|^2$ is convex for any $y\in \mathbb{R}^d$.
Then for any $Q\in {\cal P}_{ac} (\mathbb{R}^d )$, $X_i\in L^2 (\Omega, P; \mathbb{R}^d), i=1,2$, and $\lambda_1\in (0,1)$, 
\begin{equation}\label{321}
\sum_{i=1}^2\lambda_i V_S(P^{X_i},Q)-\lambda_1\lambda_2
CE[|X_1-X_2|^2]
\le V_S(P^{\sum_{i=1}^2\lambda_i X_i},Q),
\end{equation}
where $\lambda_2:=1-\lambda_1$.
Equivalently, the following is convex:
\begin{equation}\label{322}
L^2 (\Omega, P;\mathbb{R}^d)\ni X\mapsto  -V_S(P^X,Q)+CE[|X|^2]. 
\end{equation}
In particular, the following is displacement convex:
\begin{equation}\label{}
{\cal P}_{2,ac} (\mathbb{R}^d )\ni P\mapsto  -V_S(P,Q)+C\int_{\mathbb{R}^d}|x|^2P(dx).
\end{equation}
\end{theorem}

\begin{remark}\label{remark3.2}
Suppose that $a_{ij}=a_{ij} (x), \xi_i=\xi_i(x)\in C_b^\infty (\mathbb{R}^d)$
and that $a(x)$ is uniformly nonnegenerate.
Then $D_x^2\log p(0,x;1,y)$ is bounded (see \cite{3-Sheu}, Theorem B).
In particular, there exists a constant $C>0$ such that for any $y\in \mathbb{R}^d$,
$x\mapsto \log p(0,x;1,y) +C|x|^2$ is convex.
\end{remark}

For $P\in {\cal P}(\mathbb{R}^d )$,
\begin{equation}
{\cal S}(P):=
\begin{cases}
\displaystyle \int_{\mathbb{R}^d}p(x)\log p(x) dx,&P(dx)=p(x)dx,\\
\displaystyle \infty,&otherwise.
\end{cases}
\end{equation}
Let $\mu (P,Q)$ denote the joint distribution at $t=0, 1$ of the minimizer of $V_S(P,Q)$,
provided $V_S(P,Q)$ is finite.
The following implies that 
$L^2 (\Omega, P;\mathbb{R}^d)\ni X\mapsto  V_S(P^X,Q)$
 is continuous for a fixed $Q\in {\cal P}_{ac} (\mathbb{R}^d )$ such that ${\cal S}(Q)$ is finite.
The lower-semicontinuity of ${\cal P}(\mathbb{R}^d )\ni P\mapsto V_S(P,Q)$
is known and can be proved, e.g. from  Proposition \ref{pp1} and Lemma \ref{lm1}.
That of (\ref{326}) can be proved in the same way as Lemma 3.4 in \cite{1-2019}.
We give the proof for the sake of completeness.

\begin{theorem}\label{pp35}
Suppose that (A.5)-(A.6) hold. 
For $P, Q\in  {\cal P}_{2} (\mathbb{R}^d )$, 
if ${\cal S}(Q)$ is finite, then $V_S (P,Q)$ is finite and the following holds:
\begin{eqnarray}\label{325}
-V_S (P,Q)
&=&H(P\times Q| \mu (P,Q))-{\cal S}(Q)\\
&&\qquad +\int_{\mathbb{R}^d\times \mathbb{R}^d} \log p(0,x;1,y) P(dx)Q(dy).\notag
\end{eqnarray}
In particular, 
the following is weakly lower semicontinuous:
\begin{equation}\label{326}
{\cal P}_{2} (\mathbb{R}^d )\ni P\mapsto -V_S(P,Q)+C_2\int_{\mathbb{R}^d\times \mathbb{R}^d}|x-y|^2P(dx)Q(dy)
\end{equation}
(see (\ref{4.11}) for notation).
The following is also continuous  in the topology induced by $W_2$:
\begin{equation}\label{4.14}
{\cal P}_{2} (\mathbb{R}^d )\ni P\mapsto V_S(P,Q).
\end{equation}

\noindent
If ${\cal S}(Q)$ is infinite, then so is $V_S (P,Q)$.
\end{theorem}

\begin{remark}
For $C>0$ and  $P, Q\in {\cal P} (\mathbb{R}^d )$,
\begin{equation}
\Psi_{Q,C} (P):={\cal S}(P) -V_S(P,Q)+
C\int_{\mathbb{R}^d\times \mathbb{R}^d}|x-y|^2P(dx)Q(dy). 
\end{equation}
$\Psi_{Q,C} (P)$ plays a crucial role in the construction of moment measures by the SOTP
(see \cite{1-2019,Mikami} and also \cite {3-Sa} for the approach by the OTP).
Since ${\cal P}_{ac} (\mathbb{R}^d )\ni P\mapsto {\cal S}(P)$ is strictly displacement convex from Theorem 2.2 in \cite{3-McCann}, 
so is ${\cal P}_{2, ac} (\mathbb{R}^d )\ni P\mapsto \Psi_{Q,C} (P)$
under the assumption of Theorem \ref{pp34}.
\end{remark}

From Theorem \ref{pp34}, under stronger assumptions than Theorem \ref{pp35},
for a fixed $Q\in {\cal P}_{2,ac} (\mathbb{R}^d )$ such that ${\cal S}(Q)$ is finite,
we prove that 
${\cal P}_2 (\mathbb{R}^d )\ni P\mapsto V_S(P,Q)$ is Lipschitz continuous in $W_2$.

\begin{corollary}\label{co32.1}
Suppose that (A.4)' holds and 
that there exists a constant $C>0$ such that
$\log p(0,x;1,y) +C|x|^2$ is convex in $x$ for any $y\in \mathbb{R}^d$.
Then for any $Q\in {\cal P}_{2,ac} (\mathbb{R}^d )$ such that ${\cal S}(Q)$ is finite,
the following holds: 
\begin{eqnarray}\label{336AUG}
&&|V_S(P_0,Q)-V_S(P_1,Q)|\\
&\le& f (\max(||x||_{L^2(P_0)}, ||x||_{L^2(P_1)}), ||x||_{L^2(Q)} )W_2(P_0,P_1),
\quad P_0,P_1\in {\cal P}_2 (\mathbb{R}^d ),\nonumber
\end{eqnarray}
where $ ||x||_{L^2(P)}:=(\int_{\mathbb{R}^d}|x|^2P(dx))^{1/2}$, $P\in {\cal P}_2 (\mathbb{R}^d )$ and 
\begin{equation*}
f(x,y):=2C_2x^2+2(C_2y^2+C_1)+C
\end{equation*}
(see (\ref{4.11}) for notation).\\
\noindent
In particular, if $p(0,x;1,y)=(2\pi a)^{-d/2}\exp (-|y-x|^2/(2a)), a>0$, then 
\begin{eqnarray}
&&|V_S(P_0,Q)-V_S(P_1,Q)|\\
&\le &\frac{1}{2a}\{||x||_{L^2(P_0)}+||x||_{L^2(P_1)}+2 (1+\max (\sigma_0,\sigma_1))||x||_{L^2(Q)}
\}W_2(P_0,P_1),\notag
\end{eqnarray}
where 
\begin{equation*}
\sigma_i:=\biggl(\int_{\mathbb{R}^d} \biggl(x-\int_{\mathbb{R}^d}yP_i(dy)\biggr)^2P_i(dx)\biggr)^{1/2},\quad i=0,1.
\end{equation*}
\end{corollary}

We prove  Theorems \ref{pp34} and \ref{pp35}, and Corollary \ref{co32.1}.

(Proof  of Theorem \ref{pp34}).
For any $f_{i}\in C^\infty_b (\mathbb{R}^d)$, 
$u_{i} (x):=\varphi (0,x;f_{i})$ (see  (\ref{2.15}) for notation).
Then
\begin{eqnarray}\label{327}
&&\sum_{i=1}^2\lambda_i \left\{\int_{\mathbb{R}^d} f_{i} (x)Q(dx)-\int_{\mathbb{R}^d} u_i (x)P^{X_i}(dx)\right\}\\
&\le &V_S (P^{\sum_{i=1}^2\lambda_i X_i},Q)+\lambda_1\lambda_2CE[|X_1-X_2|^2].\notag
\end{eqnarray}
Indeed,
\begin{eqnarray*}
&&\sum_{i=1}^2\lambda_i \left\{\int_{\mathbb{R}^d} f_{i} (x)Q(dx)-\int_{\mathbb{R}^d} u_i (x)P^{X_i}(dx)\right\}\\
&=&\int_{\mathbb{R}^d} \sum_{i=1}^2\lambda_i f_{i} (x)Q(dx)-
E\biggl[\sum_{i=1}^2\lambda_i\{u_{i}(X_i)+C|X_i|^2\}\biggl]
+CE\biggl[\sum_{i=1}^2\lambda_i |X_i|^2\biggl],\notag
\end{eqnarray*}
\begin{eqnarray*}
&&\int_{\mathbb{R}^d} \sum_{i=1}^2\lambda_i f_{i} (x)Q(dx)\\
&\le&V_S (P^{\sum_{i=1}^2\lambda_i X_i},Q)+\int_{\mathbb{R}^d} \varphi \biggl(0,x;\sum_{i=1}^2\lambda_i f_{i}\biggl)P^{\sum_{i=1}^2\lambda_i X_i}(dx)\notag
\end{eqnarray*}
by the Duality Theorem for $V_S$ (see Corollary \ref{co2}).
\begin{eqnarray*}
&&\int_{\mathbb{R}^d} \varphi \biggl(0,x;\sum_{i=1}^2\lambda_i f_{i}\biggl)P^{\sum_{i=1}^2\lambda_i X_i}(dx)
=E\biggl[\varphi \biggl(0,\sum_{i=1}^2\lambda_i X_i;\sum_{i=1}^2\lambda_i f_{i}\biggl)\biggl]\\
&\le&E\biggl[\sum_{i=1}^2\lambda_i\biggl\{u_{i}(X_i)+C|X_i|^2\biggl\}\biggl]
-CE\biggl[\biggl|\sum_{i=1}^2\lambda_i X_i\biggl|^2\biggl].\notag
\end{eqnarray*}
In the inequality above, we considered as follows:
\begin{equation}
\label{3.46}
\varphi (t,x;f)=\log\left(\int_{\mathbb{R}^d}p(t,x;1,y)\exp (f(y)) dy\right), (t,x)\in [0,1)\times \mathbb{R}^d,
\end{equation}
\begin{eqnarray*}
&&\int_{\mathbb{R}^d} 
\exp \biggl(\log p\biggl(0,\sum_{i=1}^2\lambda_iX_i;1,y\biggl)+C\biggl|\sum_{i=1}^2\lambda_i X_i\biggl|^2+\sum_{i=1}^2\lambda_if_{i} (y)\biggl)dy\\
&\le &\int_{\mathbb{R}^d} 
\exp \biggl(\sum_{i=1}^2\lambda_i\{\log p(0,X_i;1,y)+C|X_i|^2+f_{i} (y)\}\biggl)dy\\
&\le &\prod_{i=1}^2\left(\int_{\mathbb{R}^d} 
\exp (\log p(0,X_i;1,y)+C|X_i|^2+f_{i} (y))dy\right)^{\lambda_i}
\end{eqnarray*}
by H\"older's inequality.
Taking the supremum in $f_i$ over $C^\infty_b (\mathbb{R}^d)$ on the left hand side of (\ref{327}),
the Duality Theorem for $V_S$ completes the proof (see Corollary \ref{co2}).
$\Box$

(Proof  of Theorem \ref{pp35}). 
We prove the first part.
We first prove that $V_S(P,Q)$ is finite.
Indeed, from \cite{1-RT}, 
\begin{eqnarray}\label{4.17}
V_S(P,Q)&=&
\inf\{H(\mu (dxdy)|P(dx)p(0,x;1,y)dy):\mu_1=P,\mu_2=Q\}\qquad\\
&\le& H(P(dx)Q(dy)|P(dx)p(0,x;1,y)dy)\notag\\
&=&{\cal S}(Q)-
\int_{\mathbb{R}^d\times \mathbb{R}^d}\{ \log p(0,x;1,y)\}P(dx)Q(dy)<\infty\notag
\end{eqnarray}
from (\ref{4.11}) (see (\ref{2.3}) for notation).
Here for $\mu, \nu\in {\cal P}(\mathbb{R}^d\times\mathbb{R}^d )$,
\begin{equation*}
H(\mu|\nu):=
\begin{cases}
\displaystyle \int_{\mathbb{R}^d\times \mathbb{R}^d}\left\{ \log\frac{\mu (dxdy)}{\nu (dxdy)}\right\}\mu (dxdy),&\mu\ll \nu,\\
\displaystyle \infty,&otherwise.
\end{cases}
\end{equation*}
There exists a Borel measurable $f:\mathbb{R}^d\longrightarrow \mathbb{R}$ such that
the following holds (see, e.g. \cite{1-Jamison1975}):
\begin{equation}
\mu (P,Q)(dxdy)=P(dx)p(0,x;1,y)\exp (f(y)-\varphi (0,x;f))dy
\end{equation}
(see  (\ref{3.46}) for notation).
Since $V_S(P,Q)$ is finite, $f\in L^1 (\mathbb{R}^d,P_1)$ and $\varphi (0,x;f)\in L^1 (\mathbb{R}^d,P_0)$ (see, e.g. \cite{1-RT}).
In particular, 
\begin{eqnarray}\label{4.18}
&&-V_S(P,Q)\\
&=&-\int_{\mathbb{R}^d\times \mathbb{R}^d}\left\{ \log\frac{\mu (P,Q)(dxdy)}{P(dx)p(0,x;1,y)dy}\right\}\mu (P,Q)(dxdy)\notag\\
&=&\int_{\mathbb{R}^d\times \mathbb{R}^d}(-f(y)+\varphi (0,x;f))P(dx)Q(dy)\notag\\
&=&\int_{\mathbb{R}^d\times \mathbb{R}^d}P(dx)Q(dy)
\biggl\{\log\left(\frac{P(dx)Q(dy)}{\mu (P,Q)(dxdy)}\right) -\log q(y)+\log p(0,x;1,y)\biggr\},\notag
\end{eqnarray}
which completes the proof of (\ref{325}).
$P\times Q\mapsto H(P\times Q| \mu (P,Q))$ is weakly lower semicontinuous
since $P\times Q\mapsto \mu (P,Q)$ is weakly continuous  (see \cite{1-2019})
and since $(\mu,\nu)\mapsto H(\mu| \nu)$ is weakly lower semicontinuous
(see, e.g.  \cite{1-DE}, Lemma 1.4.3).
In particular, (\ref{326}) is weakly lower semicontinuous from (\ref{325}). 
The weak lower semicontinuity of (\ref{326}) implies the upper semicontinuity of 
(\ref{4.14}) since for $P_n, P\in {\cal P}(\mathbb{R}^d), n\ge 1$,  
$W_2 (P_n,P)\to 0$ as $n\to\infty$
if and only if $P_n\to P$ weakly and $\int_{\mathbb{R}^d} |x|^2P_n(dx)\to \int_{\mathbb{R}^d} |x|^2P(dx)$
(see, e.g. \cite{1-Villani1}).
(\ref{4.14}) is also weakly lower semicontinuous by Proposition \ref{pp1} and Lemma \ref{lm1}.
We prove the last part.
\begin{equation}
q(0,x;1,y):=p(0,x;1,y)\exp (f(y)-\varphi (0,x;f)).
\end{equation}
Then by Jensen's inequality,
\begin{eqnarray}\label{4.20}
&&V_S(P,Q)\\
&=&\int_{\mathbb{R}^d\times \mathbb{R}^d}\left\{ \log q(0,x;1,y)-\log p(0,x;1,y)\right\}P(dx)q(0,x;1,y)dy\notag\\
&\ge &{\cal S}(Q)-
\int_{\mathbb{R}^d\times \mathbb{R}^d}\left\{\log p(0,x;1,y)\right\}P(dx)q(0,x;1,y)dy,\notag
\end{eqnarray}
since 
\begin{equation*}
Q(dy)=\left(\int_{\mathbb{R}^d}P(dx)q(0,x;1,y)\right)dy.
\end{equation*}
(\ref{4.11}) completes the proof.
$\Box$

\begin{remark}\label{rk11AUG}
Under (A.5)-(A.6), from Theorem \ref{thm7},
 (\ref{4.18}), and (\ref{4.20}), for $P, Q\in {\cal P}_2 (\mathbb{R}^d)$, if ${\cal S}(Q)$ is finite, then
\begin{eqnarray}
-\infty&<&-C_1+C_2^{-1}\int_{\mathbb{R}^d\times \mathbb{R}^d}|x-y|^2\mu (P,Q)(dxdy)\\
&\le&-\int_{\mathbb{R}^d\times\mathbb{R}^d}\{\log p(0,x;1,y)\}\mu (P,Q)(dxdy)\notag\\
&\le &V_S(P,Q)-{\cal S}(Q)\notag\\
&\le&-\int_{\mathbb{R}^d\times\mathbb{R}^d}\{\log p(0,x;1,y)\}P(dx)Q(dy)\notag\\
&\le & C_1+C_2\int_{\mathbb{R}^d\times \mathbb{R}^d}|x-y|^2P(dx)Q(dy)<\infty.\notag
\end{eqnarray}
\end{remark}

Remark \ref{rk11AUG} plays a crucial role in the proof of Corollary \ref{co32.1}.

(Proof of  Corollary \ref{co32.1}).
Let $X,Y\in L^2 (\Omega, P; \mathbb{R}^d)$ and $\lambda:=\min (1, ||X-Y||_2)$,
where $||X||_2:=\{E[|X|^2]\}^{1/2}$.
We prove the following when $\lambda>0$.
\begin{equation}\label{3.44AUG}
V_S(P^X,Q)-V_S(P^Y,Q)\le \lambda\{2C_2(||X||_2^2+||x||_{L^2(Q)}^2)+2C_1+C\}
\end{equation}
(see (\ref{4.11}) for notation).
From Theorem \ref{pp34},
\begin{eqnarray*}
&&(1-\lambda)V_S(P^{X},Q)+\lambda V_S(P^{\lambda^{-1}(Y-X)+X},Q)\\
&\le &\lambda(1-\lambda)C||\lambda^{-1}(Y-X)||_2^2+V_S(P^{Y},Q).
\end{eqnarray*}
since $Y=(1-\lambda )X+\lambda (\lambda^{-1}(Y-X)+X)$.
From this, 
\begin{eqnarray}\label{426rr}
&&V_S(P^X,Q)-V_S(P^Y,Q)\\
&\le &\lambda\{V_S(P^X,Q)-V_S(P^{\lambda^{-1}(Y-X)+X},Q)
+C(1-\lambda)||\lambda^{-1}(Y-X)||_2^2\}.\notag
\end{eqnarray}
Since (A.4)'  implies (A.5)-(A.6), 
\begin{eqnarray*}
&&V_S(P^X,Q)-V_S(P^{\lambda^{-1}(Y-X)+X},Q)\\
&\le&{\cal S}(Q)+C_1+2C_2(||X||_2^2+||x||_{L^2(Q)}^2)- {\cal S}(Q)+C_1
\end{eqnarray*}
from Remark \ref{rk11AUG}.
The following completes the proof of the first part:
\begin{equation*}
(1-\lambda)||\lambda^{-1}(Y-X)||_2^2
=\begin{cases}
1-\lambda, &\lambda=||X-Y||_2<1,\\
0=1-\lambda,&\lambda=1.
\end{cases}
\end{equation*}
We prove the second part.
One can set $C= (2a)^{-1}$.
From (\ref{426rr}), 
the following holds:
\begin{eqnarray}
&&V_S(P^X,Q)-V_S(P^Y,Q)\notag\\
&\le &\lambda\biggl\{ V_S(P^X,Q)-V_S(P^{\lambda^{-1}(Y-X)+X},Q)\notag\\
&&\qquad 
-\frac{1}{2a}||X||_2^2+
\frac{1}{2a}|| \lambda^{-1}(Y-X)+X||_2^2
\biggr\}+
\frac{1}{2a}(||X||_2^2-||Y||_2^2),
\notag
\end{eqnarray}
since 
\begin{equation*}
\lambda (1-\lambda)||\lambda^{-1}(Y-X)||_2^2
=\lambda(-||X||_2^2+||\lambda^{-1}(Y-X)+X||_2^2)+||X||_2^2-||Y||_2^2.
\end{equation*}
The following completes the proof: from Remark \ref{rk11AUG}, 
\begin{eqnarray*}
&&V_S(P^X,Q)-V_S(P^{\lambda^{-1}(Y-X)+X},Q)-\frac{1}{2a}||X||_2^2+
\frac{1}{2a}|| \lambda^{-1}(Y-X)+X||_2^2\\
&\le &
\frac{1}{a}\int_{\mathbb{R}^d\times\mathbb{R}^d}\langle x, y\rangle \left\{
\mu(P^{\lambda^{-1}(Y-X)+X},Q)(dxdy)-P^X(dx)Q(dy)\right\}\notag\\
&=&\frac{1}{a}\int_{\mathbb{R}^d\times\mathbb{R}^d}\langle x-E[X], y\rangle
\mu(P^{\lambda^{-1}(Y-X)+X},Q)(dxdy)\notag\\
&\le &\frac{1}{a\lambda}(||X-Y||_2 +\lambda V(X)^{1/2})||x||_{L^2(Q)}.\Box
\notag
\end{eqnarray*}


\end{document}